\documentclass[12pt]{amsart}
\usepackage{amsmath}
\usepackage{amsthm}
\usepackage{latexsym}
\usepackage{amssymb}

\newtheorem{theorem}{Theorem}

\newcounter{obsctr}

\renewcommand{\theequation}{\thesection.\arabic{equation}}
\begin{document}
\baselineskip 16pt
\def\A {{\mathcal{A}}}
\def\D {{\mathcal{D}}}
\def\R {{\mathbb{R}}}
\def\N {{\mathbb{N}}}
\def\C {{\mathbb{C}}}
\def\Z {{\mathbb{Z}}}
\def\l {\ell}
\def\ml {multline}
\def\multiline {\multline}
\def\lessim {\lesssim}
\def\d {\partial}
\def\phi{\varphi}
\def\epsilon{\varepsilon}
\title{Analytic Hypoellipticity for a Class of Sums of 
Squares of Vector Fields with Non-symplectic Characteristic
Variety}    

\author{Antonio Bove}
\address{Dipartimento di Matematica, 
Universit\`a di Bologna, Piazza
di Porta San Donato 5, 40127 Bologna, Italy}
\email{Antonio.Bove@bo.infn.it} 
\author{Makhlouf Derridj}
\address{78350 Les Loges en Josas, FRANCE}
\email{Makhlouf.Derridj@univ-rouen.fr}
\author{David S. Tartakoff}
\address{Department of Mathematics, University
of Illinois at Chicago, m/c 249, 851 S.
Morgan St., Chicago IL  60607, USA}
\email{dst@uic.edu}
\date{}
\begin{abstract}
Recently, N. Hanges proved that the operator 
$$P=\partial_t^2 + t^2\Delta_x +
 \partial_{\theta (x)}^2$$ 
in $\R^3$ is analytic hypoelliptic in the sense of germs at the origin 
and yet fails to be analytic hypoelliptic `in the strong sense' 
in any 
neighborhood of the origin (there is no neighborhood $U$ 
of the origin such that for every open subset $V$ of $U$ and 
distribution $u$ in $U,$  
$Pu$ analytic in $V$ implies that $u$ is analytic in $V$). 
Here $\partial_t=\partial / \partial{t}, \partial_{\theta(x)}=iD_{\theta(x)} = 
x_1\partial /\partial x_2 - x_2\partial /\partial x_1, 
x=(x_1,x_2),$ and $\Delta_x= \partial^2/\partial x_1^2 + 
\partial^2/\partial x_2^2.$ We give a very short $L^2$ 
proof of this result, obtained jointly with A. Bove and M.
Derridj, which generalizes easily and suggestively to other
operators.  It is striking that these operators have
non-symplectic  characteristic varieties. Finally we point
out that these results  are consistent with Treves'
conjecture. 
\end{abstract}
\maketitle
\pagestyle{myheadings}
\markboth{A. Bove, M. Derridj and D. S. Tartakoff}{Local
Analyticity for a Class of Non-symplectic Sums of Squares}
\section{Introduction and Generalizations}
\renewcommand{\theequation}{\thesection.\arabic{equation}}
\setcounter{equation}{0}
\setcounter{theorem}{0}
\setcounter{proposition}{0}  
\setcounter{lemma}{0}
\setcounter{corollary}{0} 
\setcounter{definition}{0}

In is recent paper \cite{Hanges}, Hanges considered the
operator 
\begin{equation}\label{H}
P_H=\partial_t^2 + t^2\Delta_x + \partial_{\theta (x)}^2 = \sum_1^4
X_j^2
\end{equation}
in $\R^3$ where $\partial_{\theta(x)} = 
x_1\partial /\partial x_2 - x_2\partial /\partial x_1$ and made the
interesting distinction between analytic
hypoellipticity in the germ sense and a.h.e. in the \textit{strict} 
sense. Hanges proved that the operator $ P $ was not analytic
hypoelliptic {\it in the strict sense} in any 
open set $U$ containing the origin, i.e., did not have
the property that for any open subset $V$ of $U,$ if 
$Pu$ is analytic in $V$ then so is the solution $u,$ yet 
had the property that if $Pu$ was analytic in some neighborhood
of the origin then the so was $u$ in a (possibly smaller) 
neighborhood of the origin. He ties this result to the conjecture 
of Treves concerning the Poisson strata of the operator $P,$
namely that if one writes $P=\sum_1^4 X_j^2,$ and 
considers the successive strata where 1) all $X_j$ vanish, 
2) all $X_j$ and their first brackets vanish, 3) all
$X_j$  and their first and second brackets vanish, etc.,
then the  operator should be analytic hypoelliptic in the
strict sense if and only if all these strata are symplectic.
In the case  of the particular operator being considered
here, not  even the characteristic variety is symplectic, 
being given by $t=\tau = x_1\xi_2-x_2\xi_1=0.$ 

Here we give a very elementary, and flexible, proof of the
affirmative part of his result and argue that the negative 
part is entirely reasonable as well, though we avoid entirely
the mention of so-called Treves curves, which foliate the 
characteristic variety of $P,$ in our proof. 

The generalizations we consider may be motivated by observing 
that while the ``added'' term $\partial_{\theta (x)}^2$ in $P,$ which 
suggests the celebrated non-analytic hypoelliptic example of 
Baouendi and Goulaouic, 
\begin{equation}\label{BG} P_{BG} = \partial_t^2 + t^2\partial_x^2 +
\partial_y^2 = \sum_1^3Z_j^2,
\end{equation} 
differs from this example in one essential factor---the 
integral curves of $\partial_y$ are non-compact yet those of 
$\partial_\theta$ which start close to the origin remain close. 
Hence a propagation of singularities result may be 
rephrased in terms of a germ result on analyticity. 

To put the matter differently, the $L^2$ proof of propagation
of singularities for $P_{BG}$ hinges (writing $iD=\partial$) on the
fact that in
estimating localized high derivatives of a solution $u$ in 
the $x$-direction, $\phi(x,y)D_x^p u,$ via the $L^2$ 
{\it a priori} estimate (we take $\phi$ independent of $t$ 
since for $t\neq 0,$ the operator is elliptic), 
one encounters and cannot avoid the derivation (and
bracket) 
$$\sum_j\|Z_j\phi D_x^p u\|_{L^2}^2 
\lesssim |(P_{BG} \phi D_x^p u, \phi D_x^p u)_{L^2}|
$$
$$\lesssim \sum_j\|[Z_j, \phi D_x^p] u\|_{L^2}^2  + \ldots$$
$$\lesssim \|\phi' D_x^p u\|_{L^2}^2 + \ldots $$
When the value of $j$ is $3$, i.e, we are trying to estimate 
$D_x$ derivatives of $u$ and encounter a $y$-derivative of 
$\phi$ with no gain in the number of $x$-derivatives, 
we cannot proceed, even with Ehrenpreis type 
localizing functions, to obtain analytic growth. Unless, 
of course, the $y$-derivative of the localizing function 
is supported in a region where the solution is known to be
analytic already. 

However, if the localizing function $\phi$ could be written 
as a function independent of the $y$-variable as well, this 
situation would not arise and analyticity would follow (after
some calculation, admittedly, but elementary calculations with 
no sophisticated ingredients.) 

This is what occurs when the open set under consideration 
is global in the ``$y$-direction'', as in proofs of analyticity 
which are local in some variables and global in others, as
on a tube or torus, or when the vector field $D_y$ is 
replaced by a vector field whose integral curves remain in 
any neighborhood of the point under consideration, as in 
Hanges' example, where $D_y$ is replaced by $D_\theta.$

Thus the following generalization of Hanges' example suggest 
themselves rapidly: in $(t,x)\in\R^{\ell}\times\R^k,$ and with 
$\partial_j = {\partial / \partial x_j},$
\begin{equation}\label{FC} P_1 = \Delta_t + |t|^2\Delta_x +
\sum_{i,j=1}^ka_{ij}(x,t) (x_i\partial_j-x_j\partial_i)^2
\end{equation}
for positive definite and analytic matrix valued function 
$a_{jk}.$ Note the critical feature of this operator 
that the Laplacian in $x$ commutes with each of the angular 
operators $x_i\partial_j-x_j\partial_i.$ Actually, in terms of 
estimates what is crucial is that there be a $C^\omega$ basis, 
$\{X_j\},$ of vector fields in the $x$ variables in in 
$\R^k\setminus\{0\}$ such that any bracket,
$[X_j, x_i\partial_\ell-x_\ell\partial_i]$ be a linear combination
of the angular vector fields $x_i\partial_j-x_j\partial_i$ (over 
which we have coercive control). 
 
We also remark that 
when $a_{ij}(x,t)\equiv 1,$ then it is not hard to see 
that the last sum is a constant multiple of the Laplace Beltrami operator on 
the unit sphere. 

Still more generally, let us consider $ k $ vector fields 
$ X_{1}, \ldots ,
X_{k} $ in the $x-$ variables with analytic coefficients (of $x,t$) and
$ s
$ vector fields 
$ Y_{1}, \ldots , Y_{s} $ in the $x-$ variables 
which may be singular but have analytic coefficients .  Let $
x_{0}
\in
\R^{k}
$ be a fixed point and denote by $ U \subset
\R^k$ an open 
neighborhood of $ x_{0} $. Without loss of generality we may suppose
that $ x_{0} = 0 $. We assume that

\begin{itemize}
\item[1 - ]{}
The $ \{\partial / \partial t_m, X_{j}\}_{\substack{j=1, \ldots, k\\
    m=1, \ldots,\ell}}$   
span the tangent space  on every point $ (t,x) \in
\R^{\ell}\times U$, with $x\neq 0.$
\item[2 - ]{}
$ Y_1, \ldots , Y_s $ have a compact closed family of integral
manifolds which foliate $ U $.
\item[3 - ]{}
We assume that the following commutation relations hold:
\begin{equation}
\label{1}
[\frac{\partial}{\partial t_m} , X_\ell] = \hbox{$C^\omega$ lin.
comb. of the }  Y, \frac{\partial}{\partial t}, \hbox{ and } tX, 
\;
\forall m, \ell,
\end{equation}
and for some
$C^\omega$ positive definite quadratic form $\Lambda$ in
the $X_j$'s, with coefficients independent of the $t$ variables,
\begin{equation}
\label{2}
[\Lambda, Y_{\ell}] = \hbox{$C^\omega$ quadratic
expr. in the }  Y, \frac{\partial}{\partial t}, \hbox{ and } tX,
\;
\forall \ell.
\end{equation}
\end{itemize}

Consider then the operator 
\begin{equation}
\label{GC} 
P_2 = \Delta_{t} + |t|^{2} \sum_{i,j=1}^{k} a_{ij} X_{i} X_{j} +
\sum_{i,j=1}^{s} b_{ij} Y_{i} Y_{j},
\end{equation}
where $a_{ij}(t,x)$ and $ b_{ij}(t,x) $ are $ C^{\omega} $ positive
definite matrices. We will show that we may argue as in
the particular case to obtain the result that $ P_{2} $ is
analytic hypoelliptic in the sense of germs at the origin.

Note that assumption 2 implies that we may choose a
localizing  function constant on the integral curves of $Y_1,
\ldots Y_s,$ of  Ehrenpreis type, identically
equal to one on any given  compact subset of $U$ but
vanishing outside of $U$.

We state our theorem:

\begin{theorem}
\label{th:1}
Let us consider the operator $ P_{2} $ as in (\ref{GC}), where the
coefficients $ a_{ij} $, $ b_{ij} $ are real analytic in a
neighborhood of the origin. Let $ U $ be a neighborhood of the origin
with the properties in Assumptions 1-3 above. Let $ P_{2} u = f $ hold
on the same open set $ U $, with $ f \in C^{\omega}(U) $. Then $ u $
is also in $ C^{\omega}(U) $.
\end{theorem}

\section{Proof in the case of Hanges' operator  (\ref{H})}
\renewcommand{\theequation}{\thesection.\arabic{equation}}
\setcounter{equation}{0}
\setcounter{theorem}{0}
\setcounter{proposition}{0}  
\setcounter{lemma}{0}
\setcounter{corollary}{0} 
\setcounter{definition}{0}

As remarked above, we may take localizing functions to be 
independent of $t,$ since were a derivative in $t$ to land on 
such a localizer, one would be in the region where the operator
was clearly elliptic and the analyticity 
of the solution $u$ is well known. We denote such an 
Ehrenpreis type localizing function by 
$\phi(x)=\phi_N(x)$ subject to the usual growth of its 
derivatives: $|D^\alpha \phi |\leq C^{|\alpha|+1}N^{|\alpha|}$
for $|\alpha|\leq N,$ where the constant $C$ is (universally) 
inversely proportional to
the width of the band separating the regions where $\phi \equiv
0$ and $\phi \equiv 1.$

Next, since $P$ is $C^\infty$ hypoelliptic we may assume 
that $u$ is smooth and proceed to obtain estimates for 
$D_t^pu$ and $D_{x_j}^pu$ near $0.$ 

The {\it a priori} estimate for $P,$ while subelliptic, is 
more importantly maximal: for $v\in C_0^{\infty},$ 
\setlength\multlinegap{0pt}
\begin{multline}
\|D_tv\|_{L^2}^2 + \sum_1^2\|tD_{x_j}v\|_{L^2}^2 
+ \|D_{\theta(x)}v\|_{L^2}^2\;(+ \ \|v\|_{1/2}^2) \\
\leq C|\langle Pv,v\rangle| + C\|v\|_{L^2}^2.
\end{multline}
Setting $v=\phi D_t^p u$, to begin with, we obtain 
\begin{multline*}
\|D_t\phi D_t^p u\|_{L^2}^2 + \sum_1^2\|tD_{x_j}\phi D_t^p u\|_{L^2}^2 
+ \|D_{\theta(x)}\phi D_t^p u\|_{L^2}^2\;(+ \ \|\phi D_t^p
u\|_{1/2}^2)\\[5pt] 
\leq C|\langle P\phi D_t^p u,\phi D_t^p u \rangle | + C\|\phi D_t^p
u\|_{L^2}^2 \\[5pt] 
\leq C| \langle \phi D_t^pPu, \phi D_t^pu \rangle_{L^2}| +
C\sum_1^4|([X_j^2,\phi D_t^p] u,\phi D_t^p u)|  
+ C\|\phi D_t^p u\|_{L^2}^2.
\end{multline*}
Now crucial in the brackets are the quantities (recall that we 
may take $\phi$ independent of $t,$ and clearly to localize in 
$x$ we may take it to be purely 
{\it radial} in $(x_1,x_2)$ i.e. we choose $ \phi $ to be constant on
the integral curves of $ X_{4} $), so that $X_4 \phi = 0$,
$$
[X_1,\phi D_t^p] = [X_4,\phi D_t^p] = 0,
$$
and
$$
[X_j,\phi D_t^p] = t\phi^\prime D_t^{p} - \underline{p}
\phi D_xD_t^{p-1}, \qquad j=2,3.
$$
In the first case, we may ignore the factor $t$ and recognize 
the passage from one power of $D_t$ to a derivative 
on $\phi$ as 
an acceptable swing, which, upon iteration, will lead to 
$C^{p+1}N^p \sim C^{p+1}p!$ when $p\sim N.$ The second term 
takes two powers of $D_t$ (e.g., $X_1$ from the estimate and 
one power of $D_t$ and produces a factor of $p$ and a `bad' 
vector field $D_x.$ Iterating this will yield 
$p!!D_x^{p/2}u \sim p!^{1/2}D_x^{p/2}u$ on the support of $\phi.$

On the other hand, setting $v=\phi D_{x_j}^q u,$ with perhaps 
$q=p/2,$ or, better, $v=\phi \Delta_x^{q/2} u$, where we write $
\Delta_{x} = \sum_{j} D_{x_{j}}^{2} $,
\begin{multline*}
\|D_t\phi \Delta_x^{q/2} u\|_{L^2}^2 
+ \sum_1^2\|tD_{x_j}\phi \Delta_x^{q/2} u\|_{L^2}^2 
+ \|D_{\theta(x)}\phi \Delta_x^{q/2} u\|_{L^2}^2 \\
(+ \ \|\phi \Delta_x^{q/2} u\|_{1/2}^2) \\[5pt]
\leq C|\langle P\phi \Delta_x^{q/2} u,\phi \Delta_x^{q/2} u\rangle | 
+ C\|\phi \Delta_x^{q/2} u\|_{L^2}^2 \\[5pt]
\leq C|\langle \phi \Delta_x^{q/2}Pu, \phi \Delta_x^{q/2}u \rangle_{L^2}| 
+ C\sum_1^4| \langle [X_j^2,\phi \Delta_x^{q/2}] u,\phi \Delta_x^{q/2}
u \rangle | 
\\
+ C\|\phi \Delta_x^{q/2} u\|_{L^2}^2,
\end{multline*}
and now the crucial brackets are
$$
[X_1^2,\phi \Delta_x^{q/2}] = 0, [X_4^2,\phi \Delta_x^{q/2}] = 
0
$$
and
$$
[X_j^2,\phi \Delta_x^{q/2}] = 2X_jt\phi^\prime \Delta_x^{q/2} 
- t^2 \phi^{(2)}\Delta_x^{q/2}, \qquad  j=2,3
$$ 
(where we have used rather heavily the fact that $X_4\phi =0$
since $\phi$ depends only on $x,$ and radially so, and that in 
fact $[D_\theta, \Delta_x] = 0.$) 

This last line leads to two kinds of terms, namely, for $j=2,3$,
$$
\langle 2X_jt\phi^\prime \Delta_x^{q/2}u, \phi \Delta_x^{q/2} u\rangle 
$$
and
$$
\langle t^2 \phi^{(2)}\Delta_x^{q/2} u,\phi \Delta_x^{q/2} u\rangle.
$$
Morally, these terms show the correct gain to lead to analytic
growth of derivatives, namely one must think of $t\Delta^{1/2}$
as an $X_j$ with $j=2$ or $j=3,$ and so in the first term above
one merely integrates by part noting that $X_j^*=-X_j$ and 
obtains, after a weighted Schwarz inequality, a small multiple
of the left hand side of the {\it a priori} inequality and 
the square of a term with one derivative on $\phi$ and 
$t\Delta^{1/2}$ and $q$ reduced by one, though one more 
commutator is required to make the order correct, and this 
will introduce another derivative on $\phi$ and $q$ again 
decreased by one unit, etc. The second term is of a different 
character, though the same observation reduces us essentially 
to
$$
\langle X\phi^{(2)}\Delta_x^{(q-1)/2} u,X\phi \Delta_x^{(q-1)/2}
u\rangle 
$$
in which instead of each copy of $\phi$ receiving one derivative, 
we have two derivatives on one copy and none on the other. 
Fortunately, the Ehrenpreis-type cut-off functions may be 
differentiated not merely $N$ times with the usual growth but 
$2N$ or $3N$ with no change---so in the above inner product 
we include a factor $CN$ with the copy of $\phi$ which remains
undifferentiated and a factor of $(CN)^{-1}$ with the other. 
The estimates work out just as before. 

\section{Proof in the general case (\ref{GC})}
\renewcommand{\theequation}{\thesection.\arabic{equation}}
\setcounter{equation}{0}
\setcounter{theorem}{0}
\setcounter{proposition}{0}  
\setcounter{lemma}{0}
\setcounter{corollary}{0} 
\setcounter{definition}{0}

The general case is not more complicated than the first, simplest 
case, with $\Delta_{t} = \sum_{m=1}^{\ell} D_{t_{m}}^{2}$ requiring 
us to consider each $t$-variable separately; $\Delta_x$ is replaced
by $\sum_{i,j=1}^k a_{ij}(t, x) X_{i} X_{j}$ and the square of the
angular derivative by the sum $\sum_{i, j = 1}^s b_{ij}(t, x)
Y_{i} Y_{j}$. Thus we merely give a brief sketch of the proof.

We have an \textit{a priori} estimate of the form
\begin{multline}
\label{4}
\sum_{j=1}^{\ell} \| D_{t_{j}} u \|^{2} + \sum_{j=1}^{\ell}
\sum_{h=1}^{k} \| t_{j} X_{h} u \|^{2} + \sum_{j=1}^{s} \|Y_{s} u
\|^{2} + \|u \|_{1/2}^{2} \\
\leq C \left( | \langle P_{2} u, u \rangle | + \|u \|^{2}
\right). 
\end{multline}
Let us write, as introduced above,
$$ 
\Lambda = \hbox{ a positive definite quadratic
expression in the } X_{i}
$$
and
$$ 
A = \sum_{i, j=1}^{k} a_{ij}(t, x) X_{i} X_{j}.
$$
Analogously to what has been done before we denote by $ \phi $ a
cut-off function of Ehrenpreis type, constant on the integral manifold
of the fields $Y_{1}, \ldots ,Y_{s} $ and independent of $ t $.

The problem of estimating the growth rate of the derivatives of $ u $
then reduces to estimating 
$$ 
\| \phi \Lambda^{q/2} u \|,
$$
for every natural number $ q \leq N $, where $ | \partial^{\alpha}
\phi| \leq C^{1 + |\alpha|} N^{|\alpha|} $, for $ 0 \leq |\alpha| \leq
3N$. 

We have thus to examine the structure of the commutator
$$ 
[ P_{2}, \phi \Lambda^{q/2} ] = [ \Delta_{t} + |t|^{2} A + B,
\phi \Lambda^{q/2} ],
$$
where we wrote 
$$ 
B = \sum_{i,j =1}^{s} b_{ij}(t, x) Y_{i} Y_{j}.
$$
The above quantity becomes:
\begin{multline*}
[ \Delta_{t} + |t|^{2} A + B, \phi \Lambda^{q/2} ] =  \frac{q}{2} \phi
[\Delta_{t} , \Lambda ]
\Lambda^{q/2 -1} + |t|^{2} [A, \phi ] \Lambda^{q/2} \\  + \frac{q}{2}
|t|^{2}
\phi [A,
\Lambda]
\Lambda^{q/2 -1} + \phi [B, \Lambda^{q/2}] 
= T_{1} + T_{2} + T_{3} +
T_{4},
\end{multline*}
modulo lower order terms whose treatment is easier. Let us look at
each term in the above formula, denoting by `elliptic' any term
which contains, in the inner product, two factors of the form 
maximally estimated by the operator, namely two factors each of
the form $Y, tX, \hbox{ or } \partial/\partial t.$ Such terms will be
subject to the a priori inequality (after an integration by parts) in a
recursive manner and will cause little trouble.  
\subsection{$ T_{1} $.}
Since the commutator appears in a scalar product, taking one $ t
$-derivative to the other side, we have to estimate, for some
coefficient $a(x),$ 
$$ 
q |\langle \phi [D_{t_{s}} , \Lambda ] \Lambda^{q/2 -1} u , D_{t_{s}}
\phi \Lambda^{q/2} u \rangle| \sim q |\langle \phi [D_{t_{s}} , aX ]
\Lambda^{(q-1)/2} u , D_{t_{s}}
\phi \Lambda^{q/2} u \rangle|.
$$
But by (\ref{1}), this bracket is elliptic, hence the factor of $q$
balances the decrease in the exponent of $\Lambda$ and will
iterate analytically.

\subsection{$ T_{2} $.}
We have to estimate
$$ 
|\langle |t|^{2} [\Lambda, \phi] \Lambda^{q/2} u , \phi \Lambda^{q/2}
u\rangle |.
$$
The commutator in the left hand side factor of the above scalar
product leads to an expression of the type
$$ 
2|\langle |t|^{2}  X_i\phi^\prime
\Lambda^{q/2}u, \phi \Lambda^{q/2}u\rangle|
\sim 2|\langle Z\phi^\prime
\Lambda^{(q-1)/2}u, Z\phi \Lambda^{q/2}u\rangle|
$$
with elliptic $Z$ modulo (easier) lower order terms. Here we
used just the form of $\Lambda.$ A weighted Schwarz inequality
shows that this term iterates analytically, since with the
Ehrenpreis-type localizing functions, a derivative on $\phi$
balances a decrease in $q$.
\subsection{$ T_{3} $.}
We have to estimate the scalar product
$$ 
| \langle \frac{q}{2} |t|^{2} \phi [A, \Lambda] \Lambda^{q/2 -1} u,
\phi \Lambda^{q/2} u \rangle|.
$$
Since 
\begin{multline}
[A, \Lambda] = \sum_{i,j=1}^{k} \sum_{\alpha,\beta=1}^{k} [a_{ij}
X_{i}X_{j}, aX_{\alpha}X_{\beta}] \\= \sum
\tilde{a}X^2[X,X] = \sum \tilde{a}X^2\{X \hbox{ or }
\frac{\partial}{\partial t}\}
\end{multline}
again modulo lower order terms, using assumption 3, part 1.

Now one of the $X$ factors raises $q/2-1$ to $q/2-1/2=(q-1)/2,$
with the factor of $q$ balancing the decrease from $q$ to $q-1,$
and the other two $X$'s (or one $X$ and one $\frac{\partial}{\partial
t}$ which is better, combine with $t^2$ to produce
$Z^2,$ with $Z$ elliptic. Thus this inner product, as well, iterates
analytically.

\subsection{$ T_{4} $.}
Since
$$ 
\phi [B, \Lambda^{q/2}] \sim \frac{q}{2} \phi [B, \Lambda ]
\Lambda^{q/2 - 1}, 
$$
the estimate of $ T_{4} $ boils down to computing the commutator $ [B,
\Lambda]$ and estimating the resulting terms.
But this goes as before in view of the second part of assumption 3,
since the bracket contains a product of three vector fields, two
of which are elliptic and one serves to convert $\Lambda^{q/2 -
1}$ to $\Lambda^{(q-1)/2}.$ 

This ends the proof of the Theorem in the general case.

\section{Remarks vis-\`a-vis the conjecture of Treves}
\renewcommand{\theequation}{\thesection.\arabic{equation}}
\setcounter{equation}{0}
\setcounter{theorem}{0}
\setcounter{proposition}{0}  
\setcounter{lemma}{0}
\setcounter{corollary}{0} 
\setcounter{definition}{0}

The conjecture of Treves states, in this context, that the operator 
$P$ should be analytic hypoelliptic at the origin if and only if all layers
of the Poisson stratification are symplectic. The first of these layers is 
the characteristic manifold, which is patently non-symplectic in all 
of these cases. And the operators are analytic hypoelliptic {\it in the sense
of germs} at the origin. 

The distinction is crucial. For to contradict Treves' conjecture, there 
would have to exist an open neighborhood $V$ of the origin, in which 
one could have analytic data with a non-analytic solution. And that may 
still well be the case. What we have shown is that for neighborhoods
of the origin {\it of a certain geometry relative to the operator}, 
analytic data forces analyticity of the solution.

It is well known however that in the case of the operator $ P_{2} $
there is in general propagation of the analytic wave front set (or
rather of the analytic regularity) along
the Hamilton leaves of the characteristic manifold (which are non
trivial in this case), see e.g. \cite{BoveTartakoff-tams}.

From the above proof we may see that the analyticity of the solution
is forced, in some ``adapted'' open set, by a ``global'' phenomenon,
that might be described by saying that the analytic singularities of
the solution in the open set under consideration comes from points
outside the open set lying on some Hamilton leaf of the characteristic
manifold. This is actually prevented by the ``large scale'' geometry
of the open set.

It can thus be asserted that, far from being in contradiction with
Treves' conjecture, the present result is in complete agreement with
it and points out the importance of the geometry of the non-symplectic
strata of the Poisson-Treves stratification.

\end{document}